\documentclass[final,leqno, draft]{siamltex}
\usepackage{color}

\usepackage{amsfonts}
\usepackage{amsmath}
\usepackage{amssymb}
\usepackage{graphicx}
\usepackage{bm}

\newcommand{\qed}{$\square$}

\newcommand{\curl}{{\rm curl}}
\newcommand{\dive}{{\rm div}}
\renewcommand{\grad}{{\rm grad}\,}

\newtheorem{remark}{Remark}[section]

\numberwithin{equation}{section}

\graphicspath{{../../figures/}}

\newcommand{\Nedelec}{N\'{e}d\'{e}lec }

\newcommand{\Reals}[1]{{\rm I\! R}^{#1}}
\renewcommand{\grad}{\operatorname{grad}}
\newcommand{\scalar}[2]{\langle #1, #2\rangle} 
\newcommand{\scalare}[2]{( #1 , #2 )} 
\newcommand{\scalarmu}[2]{( #1 , #2 )} 
\newcommand{\mass}[1]{\mathcal{M}_{#1}}
\newcommand{\mygrad}{\mathcal{G}}
\newcommand{\mycurl}{\mathcal{K}}
\newcommand{\myz}{\widetilde{\mathcal{Z}}}

\renewcommand{\vec}[1]{\ensuremath{\bm{#1}}}

\newcommand{\hh}{\ensuremath{\mathfrak{h}}}

\title{Numerical Approximation of Asymptotically Disappearing
  Solutions of Maxwell's Equations\thanks{Submitted to the SIAM
    Journal on Scientific Computing, 2012}}

\author{J. H. Adler\thanks{Department of Mathematics,
Tufts University, Medford, MA 02155 (james.adler@tufts.edu).} \and
V.~Petkov\thanks{Institut de Math\'{e}matiques de Bordeaux, 
351 Cours de la Lib\'{e}ration, 33405 Talence, France
(petkov@math.u-bordeaux1.fr).}
\and
L.~T.~Zikatanov\thanks{Department of Mathematics, Penn State University,
University Park, PA 16802\newline (ludmil@psu.edu).} 
}

\begin{document}
\maketitle

\begin{abstract}
  This work is on the numerical approximation of incoming solutions to
  Maxwell's equations with dissipative boundary conditions, whose
  energy decays exponentially with time. Such solutions are called
  asymptotically disappearing (ADS) and they play an important role
  in inverse back-scattering problems.  The existence of ADS is a
  difficult mathematical problem.  For the exterior of a sphere, such
  solutions have been constructed analytically by Colombini, Petkov,
  and Rauch \cite{2011ColombiniF_PetkovV_RauchJ-aa} by specifying
  appropriate initial conditions.  However, for general domains of
  practical interest (such as Lipschitz polyhedra), the existence of
  such solutions is not evident.

  This paper considers a finite-element approximation of
  Maxwell's equations in the exterior of a polyhedron, whose boundary
  approximates the sphere.  Standard N\'{e}d\'{e}lec--Raviart--Thomas
  elements are used with a Crank--Nicolson scheme to approximate the
  electric and magnetic fields.  Discrete initial conditions
  interpolating the ones chosen in
  \cite{2011ColombiniF_PetkovV_RauchJ-aa} are modified so that they
  are (weakly) divergence-free. We prove that with such initial
  conditions, the approximation to the electric field
  is weakly divergence-free for all time.  Finally, we show
  numerically that the finite-element approximations of the ADS
  also approximates this exponential decay
    (quadratically) with time when the mesh size and the time
  step become small.
\end{abstract}



\pagestyle{myheadings}%
\thispagestyle{plain}%
\markboth{ADLER, PETKOV, ZIKATANOV}{Numerical Approximation of Asymptotically Disappearing Solutions of Maxwell’s Equations}

\section{Introduction}\label{sec:intro}
This paper studies the numerical approximation of incoming solutions
to Maxwell's equations in terms of the electric field, $\vec{E}(t,\vec{x})$, and
the magnetic field, $\vec{B}(t,\vec{x})$,
\begin{eqnarray*}
&&\varepsilon \vec{E}_t - \curl\, \mu^{-1} \vec{B} = -\vec{j}, \quad
\vec{B}_t + \curl\, \vec{E} = 0,\\
&& \dive\, \varepsilon \vec{E}= 0,\quad \dive\, \vec{B}= 0,
\end{eqnarray*}
in the exterior of a spherical obstacle, with dissipative boundary
condition on the sphere (see~\eqref{eq:boundary-condition}). Here, $\varepsilon$ is the permittivity of the medium, $\mu$ is the
permeability, and $\dive\, \vec{j}=0$, where $\vec{j}$ is the known
current density of the system.  We
approximate numerical solutions, whose total energy decays
exponentially with time.  Such solutions are called {\it
  asymptotically disappearing} (ADS) and this phenomenon is of
interest for inverse back-scattering problems, since the leading term
of the back-scattering matrix becomes negligible.  Details and
construction of such solutions for the exterior of the unit sphere are
found in the recent work by Colombini, Petkov, and
Rauch~\cite{2011ColombiniF_PetkovV_RauchJ-aa}.
The dissipative boundary conditions of interest are 
\begin{equation}\label{eq:boundary-condition}
(1+\gamma )\vec{E}_{\textrm{tan}} = -\vec{n}\wedge
\mu^{-1}\vec{B}_{\textrm{tan}}\quad 
\mbox{on the boundary of the obstacle}. 
\end{equation}
Here, $\vec{n}$ is the outward unit normal to the boundary and $\gamma$ is a parameter satisfying $
\gamma > 0$. 

For a spherical obstacle, a case on which we focus on here, the boundary is
$|\vec{x}|=1$.  In~\cite{2011ColombiniF_PetkovV_RauchJ-aa}, it is
shown that for any value of the parameter $\gamma > 0$ determining the
dissipative boundary condition, there exist initial conditions such
that the boundary value problem for Maxwell's equations has a
solution, which decays exponentially in time as $\mathcal{O}(e^{rt})$,
with $r<0$. It is also interesting to note that in space such
solutions also decay asymptotically at infinity, i.e. they behave as
$\mathcal{O}(e^{r |x|})$ (see \cite{Petkov}).  Moreover, for dissipative boundary
conditions, (\ref{eq:boundary-condition}), if $\gamma > 0$, there are
no disappearing solutions, $u(T,x)$, that vanish for all $t \ge T > 0$
in the exterior of the sphere (see \cite{1986GeorgievV-aa}).  

The focus of this work is on the finite-element approximation of the
ADS in the exterior of a polyhedron that approximates the sphere.
This is a first step towards developing numerical techniques, which
later can be used to construct approximations to ADS for more
complicated obstacles and also for other symmetric hyperbolic systems
with dissipative boundary conditions. Such a general study is feasible
due to a recent result of Colombini, Petkov, and Rauch
\cite{ColombiniF_PetkovV_RauchJ-aa} for hyperbolic systems whose
solutions are described by a contraction semigroup $V(t) = e^{Gt},\: t
> 0$.  More precisely, it was shown in
\cite{ColombiniF_PetkovV_RauchJ-aa} that if certain coercivity
estimates are satisfied, then, the spectrum of the generator, $G$, in
the left half plane, $\operatorname{Re}(z) < 0$, is formed only by
discrete eigenvalues with finite multiplicities. Every such
eigenvalue, $\lambda\,: Re(\lambda) < 0$, yields an ADS solution,
$u(t,x) = e^{\lambda t}f(x)$, with $G f = \lambda f$.  We notice that
in \cite{2011ColombiniF_PetkovV_RauchJ-aa} only real eigenvalues of
the generator $G$ have been constructed. The question of the existence
of complex eigenvalues is open and numerical results can shed light on
the open question of existence of complex eigenvalues.

On the other hand, the {\it existence} and the {\it location} of
eigenvalues of $G$ for less regular obstacles is a difficult and
interesting mathematical problem (both from theoretical and numerical
points of view), albeit it falls beyond the scope of the work reported
here.

For the finite-element spaces that we use, their properties and
implementation in the numerical models based on Maxwell's system are
more or less known. Classical references on the piecewise polynomial
spaces relevant in such approximations are the papers by Raviart and
Thomas~\cite{1977RaviartP_ThomasJ-aa} (for two spatial dimensions),
\Nedelec \cite{1980NedelecJ-aa,1986NedelecJ-aa}, and
Bossavit~\cite{1988BossavitA-aa}. The method that is used here is an
application of the techniques developed by
Brezzi~\cite{1974BrezziF-aa}, (see also Brezzi and
Fortin~\cite{1991BrezziF_FortinM-aa}). Many results and references on
Maxwell's system and its numerical approximation are found in
Hiptmair's work~\cite{2002HiptmairR-aa} and in a monograph by
Monk~\cite{2003MonkP-aa}. In most of these works, the emphasis is on
systems with perfect conductor boundary conditions.  
The important difference between our work and the previous ones is
that, here, we apply the methods to a problem with dissipative boundary conditions, where the
electric field, $\vec{E}$, and the magnetic field, $\vec{B}$, are paired
together.  A related discussion on the role of the boundary conditions
in the numerical solution of systems of PDEs are discussed in a recent
paper by~Arnold, Falk and Gopalakrishnan~\cite{AFG11}.  

We consider the finite-element problem associated with Maxwell's
equations in a finite domain between two spheres,
$\Omega=\{\vec{x}\ |\ 1< |\vec{x}|< R\}$, for a fixed, large enough
$R$.  The discretization of Maxwell's equations
that keeps the coupling between the electric and the magnetic
field intact, and can be derived via the modern techniques in finite
element exterior calculus, is described in Arnold, Falk and
Winther~\cite{Arnold1,2010ArnoldD_FalkR_WintherR-aa}.  
We obtain a
computational domain (denoted again with $\Omega$), which is
\emph{polyhedron} and which \emph{approximates} the exterior of
the annular domain $\{\vec{x}\ |\ 1< |\vec{x}|< R\}$. Note that the ADS
constructed in~\cite{2011ColombiniF_PetkovV_RauchJ-aa} for the
exterior of the sphere do not need to satisfy the dissipative boundary
condition on the polyhedron. However, we show numerically that the
finite-element solution with initial conditions approximating those in
\cite{2011ColombiniF_PetkovV_RauchJ-aa} yields good approximation of
the ADS constructed analytically, when the mesh size becomes small and
the polyhedron gets closer to the sphere.


This paper is organized as follows.  In Section \ref{sec:background},
we introduce some notation and state the strong form of the boundary
value problem for Maxwell's equations that we consider.  Section
\ref{sec:weakform} describes the variational (weak) formulation and
discusses the energy decay of the corresponding system.  Next, we
discuss the discretization of this variational form and how we can
guarantee a good approximation of the ADS in Section
\ref{sec:discrete}.  In Section \ref{sec:time}, we describe the matrix
representation of the semi-discrete system and the properties of the
Crank--Nicolson scheme that we use for time stepping.  Numerical
results for the sphere are shown in Section \ref{sec:numerics}.
Finally, concluding remarks and discussions on
constructing initial conditions as well as a choice of parameters for
more complicated obstacles are presented in Section
\ref{sec:conclusion}.

\section{Notation and preliminaries}\label{sec:background}
First, some standard notation is introduced, which is needed in the following
sections. The Euclidean scalar product between two vectors in $\Reals{d}$  is denoted by
$\scalar{\vec{a}}{\vec{b} }$ and the corresponding norm is   $|\vec{a}|^2=\scalar{\vec{a}}{\vec{a}}$.
The standard $L^2(\Omega)$ scalar product and norm are denoted by
$(\cdot,\cdot)$ and $\|\cdot\|$,  respectively. 
\subsection{Maxwell's system}\label{subsect:model}
The system of partial differential equations (PDEs) of interest is
Maxwell's system with a dissipative boundary condition (impedance
boundary condition).  Let $\mathcal{O}$ be a bounded, connected (could
be convex) domain, $\mathcal{O}\subset \Reals{3}$. As stated
in the introduction, we consider Maxwell's equations in the exterior of $\overline{\mathcal{O}}$, that is, in
$\Omega=\Reals{3}\setminus \overline{\mathcal{O}}$, which is as follows:
\begin{eqnarray}
\varepsilon \vec{E}_t - \curl\, \mu^{-1} \vec{B} &=& -\vec{j}, \label{eq:EcB0}\\ 
\vec{B}_t + \curl\, \vec{E} &=& 0,  \label{eq:BcE0}\\
\dive\, \varepsilon \vec{E}& =& 0\label{eq:divD0},\\
\dive\, \vec{B}& =& 0 \label{eq:divB0}.
\end{eqnarray}

For the rest of the paper, we assume that the permitivity,
$\varepsilon$, and the permeability, $\mu$, are equal to $1$.  Future work will involve
investigating numerical methods when these parameters are allowed to
vary with the domain.  We also assume that $\Omega$ is bounded, which means
$\Omega=\mathcal{S}\setminus \overline{\mathcal{O}}$,
where $\mathcal{S}$ is a ball in $\Reals{3}$ with sufficiently large
radius. While in general this could be a restriction, in this case it
is not, since the solutions that are approximated decay exponentially
when $|\vec{x}|\to\infty$.  

\subsection{Dissipative boundary conditions\label{sec:boundary-conditions}}

The boundary conditions that are of interest are also known as impedance boundary conditions.  To state such
type of boundary conditions, we first define 
$$
\Gamma_{i}=\partial\Omega\cap \partial \mathcal{O},\qquad
\Gamma_{o}=\partial\Omega\setminus \Gamma_i,
$$
where $\Gamma_i$ represents the boundary of the inner obstacle, and
$\Gamma_o$ the outer boundary (i.e., the boundary of $\mathcal{S}$).
The orthogonal projection, $Q_\textrm{tan}$, on the component of a
vector field tangential to $\Gamma_{i}$ is also needed, which for any
$\vec{x}\in\Gamma_{i}$ and a vector field $\vec{F}(\vec{x})\in \Reals{3}$ is defined as the
tangential component of $\vec{F}(\vec{x})$, namely:
$$
\vec{F}_\textrm{tan}=Q_\textrm{tan} \vec{F} = \vec{F}-\langle \vec{F},\vec{n}\rangle \vec{n}=-\vec{n}\wedge\left(\vec{n}\wedge \vec{F}\right),
$$
where $\vec{n}$ is the normal vector to the surface $\Gamma_i$.  We
denote by $\vec{n}$ the normal pointing away from the domain
(i.e. pointing into $\mathcal{O}$).
We note that all the quantities above depend on $\vec{x}\in \Gamma_{i}$.
The boundary condition of interest is the one
in~\eqref{eq:boundary-condition} and it is recalled here:
\begin{equation*}
(1+\gamma )\vec{E}_{\textrm{tan}} = -\vec{n}\wedge\vec{B}_{\textrm{tan}}\quad \mbox{or equivalently}\quad 
(1+\gamma )\vec{E}_{\textrm{tan}}=  -\vec{n}\wedge\vec{B}.
\end{equation*}
As pointed out above, $\gamma > 0$  is a constant, i.e.  $\gamma\in
\Reals{}$.  However, the same methods can be applied
for $\gamma(x) > 0$ as a function on the boundary of the domain.  

\begin{remark} \label{rem:perfect}
Note that for a perfectly conducting obstacle, the tangential component of $\vec{E}$ vanishes
on the boundary, namely:
\begin{equation*}
\vec{E}\wedge \vec{n}  =  \vec{0}, \quad \vec{x}\in \partial\Omega.
\end{equation*}
However, again, this paper considers the case of an
impedance condition, where the obstacle is not a perfect conductor.
This is closer to real-world applications, where
dissipative boundary conditions occur frequently. 
\end{remark}


\section{Function spaces and variational formulation}\label{sec:weakform}
\subsection{Function spaces} To approximate the differential
problem~\eqref{eq:EcB0}--\eqref{eq:divB0} with the boundary conditions
given in~\eqref{eq:boundary-condition},
the function spaces for the problem at hand need to be identified.  Given a Lipschitz
domain, $\Omega$, and a differential operator, $\mathcal D$, a
standard notation for the following spaces is used:
$$ 
H({\mathcal D})= \{\vec{v}\in (L^2(\Omega))^{d_1}, {\mathcal D} \vec{v}\in
 (L^2(\Omega))^{d_2}\},
$$ 
with the associated graph norm, 
$$
\|\vec{u}\|_{\mathcal D; \Omega}^2 = \|\vec{u}\|^2 + \|\mathcal D\vec{u}\|^2,
$$
where $d_1$ is the dimension of the problem (three for this paper) and
$d_2$ depends on the operator.  By taking $\mathcal D=\dive$ ($d_2 =
1$) or $\mathcal D=\curl$ ($d_2 = 3$), the Sobolev spaces
$H(\dive)$ and $ H(\curl)$ are obtained.  Also, notice that 
$$
H^1(\Omega)=H(\grad), \qquad L^2(\Omega)=H(id) \qquad d_2 = 3.  
$$
For example,  $H(\curl)$ is the space of $L^2(\Omega)$ vector-valued functions,
whose $\curl$ is also in $L^2(\Omega)$. Similarly for $H(\grad)$
and $H(\dive)$. 

The following three spaces are needed (the first one for scalar
functions and the second and third for vector-valued functions): 
\begin{eqnarray*}
  H_0(\grad) = H^1_0(\Omega) & = & \{ v\in H^1(\Omega)\quad\mbox{such
    that}\quad v\big|_{\partial\Omega}=0\},\\
  \widetilde{H}_{\textrm{imp}}(\curl) & = & \{\vec{v}\in H(\curl)\quad\mbox{such
    that}\quad \vec{v}\wedge \vec{n}\big|_{\Gamma_o}=0\},\\
 H_{\textrm{imp}}(\dive) & = & \{ \vec{v}\in H(\dive)\quad\mbox{such that}\quad \scalar{\vec{v}}{\vec{n}}\big|_{\Gamma_o}=0\}.
\end{eqnarray*}
Note that the tangential component of the electric
  field, $\vec{E}$, on 
  $\Gamma_o=\partial\Omega\setminus \Gamma_{i}$, is set to zero for
  the elements of $\widetilde{H}_{\textrm{imp}}(\curl)$.  From this
  and the fact that $\vec{B}_t = -\curl~\vec{E}$, it follows that the
  normal component of the magnetic field $\vec{B}$ is also zero on
  this outer boundary, which is the boundary condition incorporated in
  $H_{\textrm{imp}}(\dive)$.  

Another issue to address is related to the fact that the boundary of
the computational domain consists of two connected components
$\Gamma_i$ and $\Gamma_o$, which is an artifact of how we treat the
far field.  Even though there are methods, such as perfectly matched
layer, absorbing boundary conditions, or others, that treat the far field
more accurately, the solutions presented here decay rapidly in space
and, therefore, we can simply add an outer boundary with homogeneous
Dirichlet boundary conditions.  
If the initial conditions is a harmonic form for $\vec{E}$ and
$\vec{B}_0=0$, then, as is easily computed, the solution does not change in
time. If we did not have the outer boundary, then such phenomenon will
not happen. From this, it can be concluded that if a harmonic form is
part of the initial condition, then this part is going to be unchanged 
when propagated in time.  However, since the goal is to show that the 
energy of the system dissipates over time (for any initial condition
that does not have the harmonic form as a component), we consider only
initial conditions that are orthogonal to the space of harmonic forms.  To
resolve the ambiguity, we consider electric fields in a subspace of
$\widetilde{H}_{\textrm{imp}}(\curl)$, which is orthogonal to the
one-dimensional space of harmonic forms.  Let $\hh$ be the unique
solution to the Laplace equation:
\begin{equation}\label{eq:harmonic-form}
-\Delta \hh=0,\quad \hh=1 \quad\mbox{on}\quad\Gamma_i,\quad 
\mbox{and}\quad \hh=0 \quad\mbox{on}\quad\Gamma_o.
\end{equation}
Then, define $H_{\textrm{imp}}(\curl)$ as the space of functions orthogonal to $\grad\hh$:
\begin{equation}\label{eq:himp-definition}
  H_{\textrm{imp}}(\curl) = 
\{\vec{v}\in   \widetilde{H}_{\textrm{imp}}(\curl) \quad\mbox{such that}\quad (\vec{v},\grad\hh)=0\}.
\end{equation}
Finally, for the time-dependent problem considered here, the relevant function spaces
are 
\begin{eqnarray*}
H_0(\grad;t) & = & \{ v(t,\cdot)\in H_0^1(\Omega) \quad\mbox{for all $t\ge 0$}\},\\
H_{\textrm{imp}}(\curl;t) & = & \{ \vec{v}(t,\cdot)\in
H_{\textrm{imp}}(\curl), \quad\mbox{for all $t\ge 0$}\},\\
H_{\textrm{imp}}(\dive;t) & = & \{ \vec{v}(t,\cdot)\in H_{\textrm{imp}}(\dive), \quad\mbox{for all $t\ge 0$}\}.
\end{eqnarray*}
In another words, if $H(\mathcal{D})$ denotes any of the Hilbert
spaces $H_0(\grad)$, $H_{\textrm{imp}}(\curl)$, or $H_{\textrm{imp}}(\dive)$, then
$H(\mathcal{D};t)$ is the space of functions, which for each $t\in
[0,\infty)$ takes on values in $H(\mathcal{D})$. We assume that the
elements of any of the spaces $H_0(\grad;t)$,
(resp. $H_{\textrm{imp}}(\curl;t)$, or $H_{\textrm{imp}}(\dive;t)$) are
differentiable with respect to $t$ as many times as needed.  We refer
to McLean~\cite{2000McLeanW-aa} and also Monk~\cite{2003MonkP-aa} for properties of the above spaces, and
related density results.


\subsection{Variational formulation} 
Next, we derive a weak form which was shown to us by
D.~N.~Arnold~\cite{2011ArnoldDN-aa}.  
We introduce $p\in H_0(\grad;t)$ such that 
\begin{equation}\label{eq:pre-1}
  (p_t,q) = \scalare{\vec{E}}{\grad q},\quad\mbox{for all}\quad q\in
  H_0(\grad), \quad p(0,\vec{x})=0.
\end{equation}
This is an auxiliary variable associated with the constraint that
$\vec{E}$ is divergence-free.  If the initial condition for $\vec{E}$
is divergence-free, then as shown below, $p$ is zero for all
times. However, if the initial guess is not divergence-free, then $p$
does not have to be zero. 

First, multiply equations~\eqref{eq:BcE0} and \eqref{eq:EcB0} 
by test functions $\vec{C}\in H_{\textrm{imp}}(\dive)$ 
and $\vec{F}\in H_{\textrm{imp}}(\curl)$, respectively.  Next, integrate by parts and use
boundary conditions~\eqref{eq:boundary-condition} and the
identities
\[
\scalar{\vec{n}\wedge\vec{B}}{\vec{F}}
=\scalar{(\vec{n}\wedge
  \vec{B}_\textrm{tan})}{\vec{F}_\textrm{tan}},\quad 
\scalar{\vec{E}_{\textrm{tan}}}{\vec{F}_{\textrm{tan}}}=\scalar{\vec{n}
  \wedge \vec{E}}{\vec{n}\wedge \vec{F}},
\]
to obtain for all $\vec{C}\in H_{\textrm{imp}}(\dive)$, and for all
$\vec{F}\in H_{\textrm{imp}}(\curl)$,
\begin{equation}\label{eq:pre-2}
\begin{array}{rcl}
&&\scalarmu{\vec{B}_t}{\vec{C}} + \scalarmu{\curl\, \vec{E}}{\vec{C}} = 0,\\
&&\scalare{\vec{E}_t}{\vec{F}} + \scalare{\nabla p}{\vec{F}}-\scalarmu{\vec{B}}{\curl\, \vec{F}} \\
&&~~~~~~~~ + (1+\gamma)\int_{\Gamma_{i}}\scalar{\vec{n} \wedge
  \vec{E}}{\vec{n}\wedge \vec{F}}\;d\gamma = -(\vec{j},\vec{F}), 
\end{array}
\end{equation}
Finally, we get the following variational problem:\\
Find $(\vec{E},\vec{B},p)
\in H_\textrm{imp}(\curl;t)\times H_{\textrm{imp}}(\dive;t)\times H_0(\grad;t)$,
such that for all 
$(\vec{F},\vec{C},q) \in H_\textrm{imp}(\curl)\times
H_{\textrm{imp}}(\dive)\times H_0^1(\Omega)$ and for all $t>0$,
\begin{eqnarray}
&&\scalare{\vec{E}_t}{\vec{F}} =-\scalare{\grad  p}{\vec{F}}
+\scalarmu{\vec{B}}{\curl\, \vec{F}} - (1+\gamma)\int_{\Gamma_{i}}\scalar{\vec{E}_\textrm{tan}}{\vec{F}_\textrm{tan}}
-(\vec{j},\vec{F}),\label{mixed-form-1}\\
&&\scalarmu{\vec{B}_t}{\vec{C}} =-\scalarmu{\curl\, \vec{E}}{\vec{C}},\label{mixed-form-2}\\
&&(p_t, q)= \scalare{\vec{E}}{\grad q}.\label{mixed-form-3}
\end{eqnarray}
At $t=0$, the following initial conditions are needed,
\begin{equation}\label{eq:initial-condition}
\vec{E}(0,\vec{x})=\vec{E}_0(\vec{x}),\quad
\vec{B}(0,\vec{x})=\vec{B}_0(\vec{x}),\quad
 p(0,\vec{x})=0.
\end{equation}
Here, $\vec{E}_0\in H_{\textrm{imp}}(\curl)$, $\vec{B}_0\in
H_{\textrm{imp}} (\dive)$, and we assume that $\vec{B}_0$ is
divergence-free.  As it is well known, this assumption implies that
$\vec{B}$ is divergence-free for all $t>0$.  Even
  though there are no derivatives acting on $\vec{B}$ in the
  variational form, we choose $\vec{B} \in H_{\textrm{imp}}(\dive)$,
  because $\dive~\vec{B}=0$ is enforced strongly and, therefore,
  $\vec{B} \in H(\dive)$.  

\begin{remark}
  Regarding the existence and uniqueness of the solution to the
  variational problem,~\eqref{mixed-form-1}--\eqref{mixed-form-3}, we
  point out that in the case of a perfect conductor, solution results for
  Lipschitz polyhedra are found in Pauly and
  Rossi~\cite{2011PaulyD_RossiT-aa} and also in Birman and
  Solomyak~\cite{1987BirmanM_SolomyakM-aa,1989BirmanM_SolomyakM-aa}.
  However, for the dissipative boundary conditions we consider and for
  Lipschitz polyhedral domains, to our best knowledge, there are no
  results on existence and uniqueness. For solution results on smooth
  domains ($C^1$) we refer to the monograph by
  Petkov~\cite{1989PetkovV-aa}, where the general case of symmetric
  hyperbolic systems is considered and existence and uniqueness
  are derived using semigroup theory. For Lipschitz polyhedral domains
  it is not straightforward and difficult to see that the solution is
  related to a contraction semigroup. 
\end{remark}

Next, the following proposition shows that for a divergence-free
$\vec{E}_0\in H_{\textrm{imp}}(\curl;t)$, the variational form
satisfies the divergence-free condition for $\vec{E}$ weakly for all
time.  
\begin{proposition}\label{prop:weakly-div-free} Let
  $u=(\vec{E},\vec{B},p)$ satisfy
  Equations~\eqref{mixed-form-1}--\eqref{mixed-form-3} and the initial
  conditions~\eqref{eq:initial-condition}. If $\vec{E}_0\in
  H_{\textrm{imp}}(\curl)$ is weakly divergence-free, then, for all $t$ and $\vec{x}$,
  $p(t,\vec{x})=0$ and, hence,
\[
\scalare{\vec{E}}{\grad      q}=0,\quad\mbox{for all}\quad q\in H_0^1(\Omega).
\]
\end{proposition} 
{\bf Proof.~~} In Equation~\eqref{mixed-form-3}, take
  $q\in H_0^1(\Omega)$, then differentiate with respect to $t$ to get
\begin{equation}\label{eq1}
(p_{tt}, q)  = \scalare{\vec{E}_t}{\grad q},
\quad\mbox{for all}\quad q\in H_0^1(\Omega).
\end{equation}
In the first equation, \eqref{mixed-form-1}, set $\vec{F}=\grad
q$.  Note that the tangential component of $\grad q$ is zero, because
$q\in H^1_0(\Omega)$ has a vanishing trace on the boundary of
$\Omega$. Thus, the boundary term vanishes and we are
  left with
$$(\vec{E}_t,\grad q) = -(\grad p,\grad q) + (\vec{B},\curl\grad q) -
(\vec{j},\grad q).$$
Using integration by parts on the last term yields
$$(\vec{j},\grad q) = \int_{\partial\Omega}(\vec{n}\cdot\vec{j})q
-(\dive\,\vec{j},q) = 0,$$
since $\vec{j}$ is divergence-free and $q$ vanishes on the
boundaries.  Taking this into account, along with~\eqref{eq1} and the identity
$\curl\,\grad=0$,
rewrite Equation~\eqref{mixed-form-1} as follows:
\begin{equation}\label{eq2}
(p_{tt},q) +\scalare{\grad p}{\grad q}=0,
\quad\mbox{for all}\quad q\in H_0^1(\Omega).
\end{equation}
Since $p(0)=0$ and $p_t(0)=-\dive\,\vec{E}(0)=0$, it is concluded
that $p(t,x)=0$ for all $t>0$ and all $\vec{x}\in \Omega$, since it is
a solution of the homogenous wave equation, i.e., Equation
~\eqref{eq2}.  The condition $p_t(0) = -\dive\,\vec{E}(0)$ comes from
  integrating \eqref{mixed-form-3} by parts and noting that $q$ vanishes
  on the boundary.  Then, from \eqref{mixed-form-3}, the desired result for $\vec{E}$ follows.\hfill\qed

\noindent Next, we show an energy estimate. 
\begin{proposition}\label{the-only-proposition} Let
  $u=(\vec{E},\vec{B},p)$ satisfy
  Equations~\eqref{mixed-form-1}--\eqref{mixed-form-3} and the initial
  conditions~\eqref{eq:initial-condition}. Assume that
  $\vec{j}=\vec{0}$ (i.e., no external forces). Then, the following
  estimate holds for all $T\geq 0$:
\[
\|p(T,\cdot)\|^2+\|\vec{E}(T,\cdot)\|^2 +\|\vec{B}(T,\cdot)\|^2
\le \|\vec{E}_0\|^2 +\|\vec{B}_0\|^2
\]
\end{proposition} 
{\bf Proof.~~} Fix $t$ and take
$(q,\vec{F},\vec{C})=(p,\vec{E},\vec{B})$. Summing up the three
equations~\eqref{mixed-form-1}--\eqref{mixed-form-3} gives
\begin{equation}
\frac{d}{dt} 
\left(\|p\|^2+\|\vec{E}\|^2 +\|\vec{B}\|^2\right)= 
-2(1+\gamma)\|\vec{E}_\textrm{tan}\|^2_{L^2(\Gamma_{i})}
\end{equation}
This identity holds for any $t$ and the proof is concluded  after integrating with respect to
time. \hfill\qed


\section{Finite-element discretization}\label{sec:discrete}
\subsection{Domain partitioning} To devise a discretization
of~\eqref{mixed-form-1}--\eqref{mixed-form-3}, 
we approximate
the exterior of the sphere by a polyhedral domain, which is decomposed as a union of
simplices (tetrahedrons). The polyhedral domain and its splitting is
obtained by mapping a corresponding splitting of a cube to a
polyhedron with vertices on the sphere.

\begin{figure}[h!]
\begin{minipage}[b]{0.4\linewidth}
\centering
\includegraphics[scale=0.15]{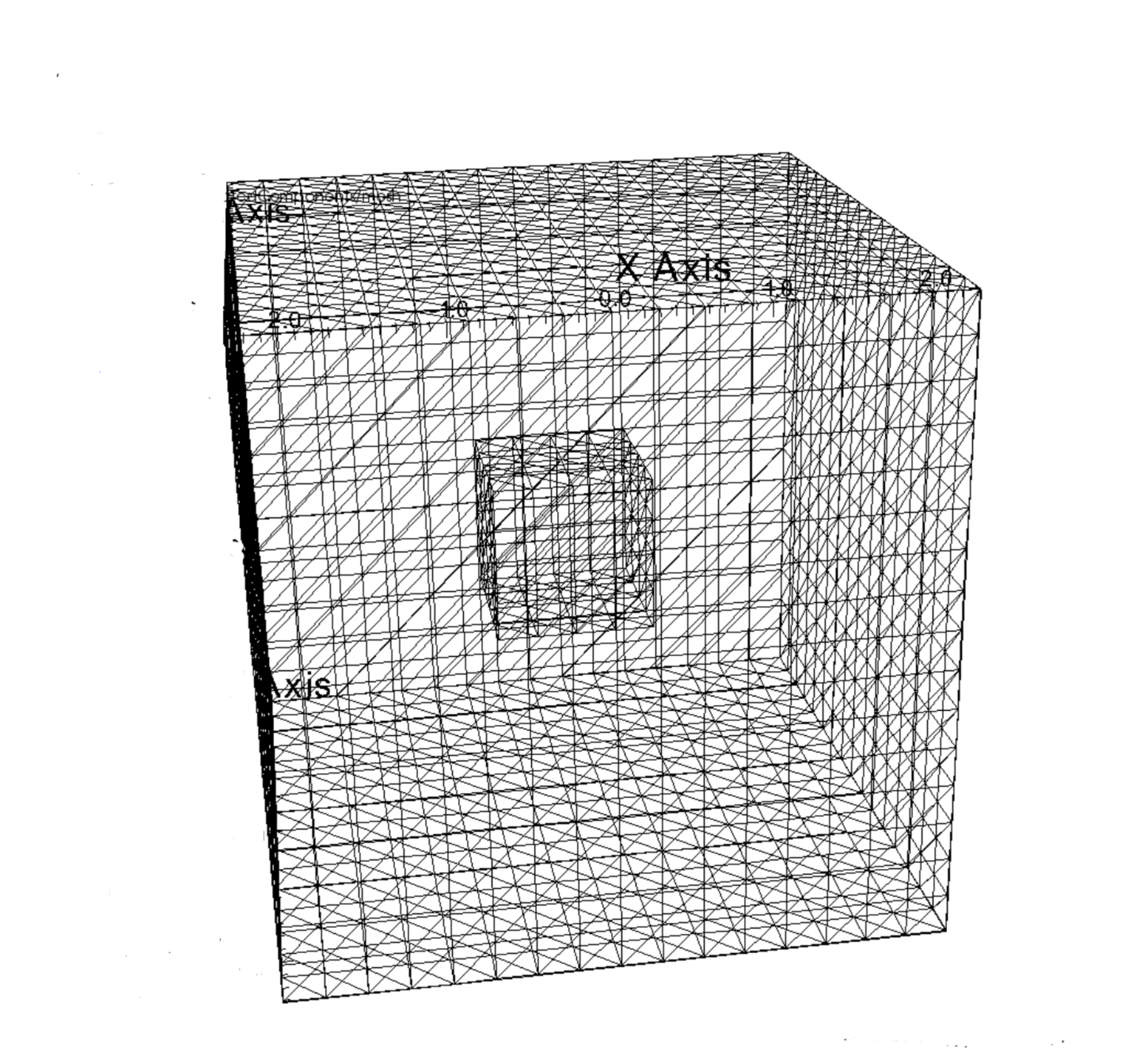}
\caption{Cube Obstacle}
\label{fig:cube}
\end{minipage}
\hfill
\begin{minipage}[b]{0.4\linewidth}
\centering
\includegraphics[scale=0.21]{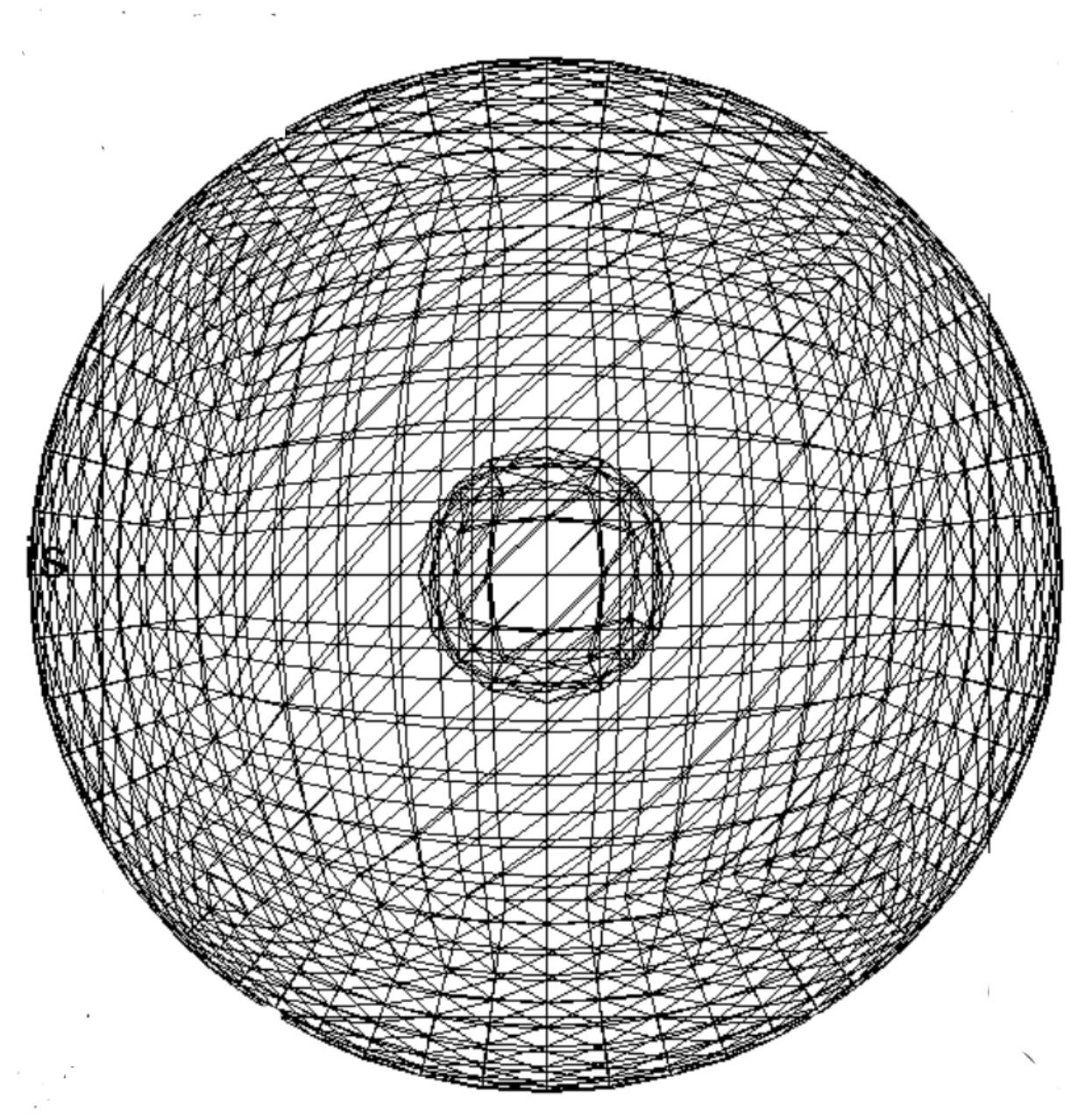}
\caption{Sphere Obstacle}
\label{fig:sphere}
\end{minipage}
\end{figure}

We consider a cube $\widetilde{\Omega}=(-R/2,R/2)^3$ and we split it
into $\frac{R^3}{h^3}$ cubes each with side length $h$ and integer $R>1$. Here $h=2^{-J}$,
for some $J\geq 2$. From this partition, we remove all cubes that have
nonempty intersection with the open cube
$\widetilde{\omega}=(-1/2,1/2)^3$.  Finally, we split
each of the cubes from the lattice into 6 tetrahedrons and this gives
a partitioning of
$\overline{\widetilde{\Omega}\backslash\widetilde{\omega}}$ into
simplices. Note that this partition has the same vertices as the
lattice and
\[
\overline{\widetilde{\Omega}\backslash\widetilde{\omega}}
=\cup_{K\in \mathcal{T}_h} \overline{K}. 
\]
From this, we obtain a polyhedron approximating the region between two
spheres by mapping a vertex of the lattice whose Cartesian coordinates
are $\vec{x}$ and whose spherical coordinates are
$(|\vec{x}|_{\ell_2},\theta,\phi)$ to the point with spherical
coordinates $(|\vec{x}|_{\ell_\infty},\theta,\phi)$.
Clearly this maps the
interior boundary of
$\overline{\widetilde{\Omega}\backslash\widetilde{\omega}}$ to the
unit sphere, and the outer boundary to a sphere with radius
$R$. An example is shown in Figures~\ref{fig:cube}--\ref{fig:sphere}. We
note that when $h\to 0$ the corresponding polyhedron approximates the
region between the
spheres. 

\subsection{Finite-element spaces}\label{subsect:FE}
For the discrete problem, we use standard piecewise linear continuous
elements paired with \Nedelec \cite{1980NedelecJ-aa,1986NedelecJ-aa}
finite-element spaces.  We introduce the spaces
$H_{h,\textrm{imp}}(\curl)\subset
\widetilde{H}_{\textrm{imp}}(\curl)$, and $H_{h,{\textrm{imp}}}(\dive)\subset
H_{\textrm{imp}}(\dive)$ and the FE solution is denoted by
$(\vec{E}^h,\vec{B}^h,p^h)\in H_{h,\textrm{imp}}(\curl)\times
H_{h,{\textrm{imp}}}(\dive)\times H_h(\grad)$.  More precisely, we define
\begin{eqnarray*}
H_{h,{\textrm{imp}}}(\dive)=\{\vec{C}\in H_{h,{\textrm{imp}}}(\dive),\ 
\vec{C}|_K=\vec{a}_K+\beta_K \vec{x}, \ \forall K\in \mathcal{T}_h\},\\
H_{h,0}(\grad)=\{q\in H_0^1(\Omega),\ 
q|_K\ \mbox{is linear in $\vec{x}$}, \ \forall K\in \mathcal{T}_h\}.
\end{eqnarray*}
The space $H_{h,{\textrm{imp}}}(\curl)$ is a properly chosen subspace
of $\widetilde{H}_{h,\textrm{imp}}(\curl)$,
which is orthogonal to the gradients of the \emph{discrete} harmonic
forms (but not necessarily to $\grad \hh$).  This
intermediate space is defined as
\[
\widetilde{H}_{h,{\textrm{imp}}}(\curl)=\{\vec{F}\in \widetilde{H}_{{\textrm{imp}}}(\curl),\ 
\vec{F}|_K=\vec{a}_K + \vec{b}_K\wedge  \vec{x}, \ \forall K\in
\mathcal{T}_h\}.
\]
Next, the discrete harmonic form, $\hh^h$, is defined as the
unique piecewise linear continuous function satisfying
\begin{eqnarray*}
(\grad \hh^h,\grad q) &=& 0,
\quad \mbox{for all}\quad q\in H_{h,0}(\grad),\\
\hh^h&=&1 \quad\mbox{on}\quad\Gamma_i,\\
\hh^h&=&0 \quad\mbox{on}\quad\Gamma_o. \\
\end{eqnarray*}
Then,
\[
H_{h,\textrm{imp}}(\curl)=\{\vec{F}\in \widetilde{H}_{h,\textrm{imp}}(\curl),
\ (\vec{F},\grad\hh^h)=0\}.
\]
In the definitions above, $\vec{a}_K, \vec{b}_K\in \Reals{3}$ are
constant vectors for every simplex $K$ in the partition and
$\beta_K\in \Reals{}$. 

Spaces corresponding to the time-dependent problem are analogously
defined using the definitions from the previous section.  We denote
these spaces by $H_{h,\textrm{imp}}(\curl;t)$, $H_{h,{\textrm{imp}}}(\dive;t)$ and
$H_{h,0}(\grad;t)$, respectively.  Note that $\vec{v}\in
H_{{\textrm{imp}}}(\curl;t)$ implies that the tangential components of
$\vec{v}$ are continuous.  The other spaces induce certain
compatibility conditions as well.  For example, the requirement $\vec{B}\in
H(\dive;t)$ in the definition of the space $H_{h,{\textrm{imp}}}(\dive;t)$ at the
beginning of Section~\ref{subsect:FE} is equivalent to saying that the
\emph{normal} components of the elements from $H_{h,{\textrm{imp}}}(\dive;t)$ are
continuous across element faces. It is also easy to check that $q\in
H_{h,0}(\grad;t)\subset H_0^1(\Omega)$, in the definition above, implies
that $q$ is a continuous function (because it is a piecewise polynomial
function, which is in $H_0^1(\Omega)$). 

Finally, it is important to note that
Equations~\eqref{mixed-form-1}--\eqref{mixed-form-3} make sense for
all $\vec{E}\in H(\curl)$. As stated above, the piecewise polynomial
functions on tetrahedral partitions of $\Omega$ are in $H(\curl)$ if
their tangential components across the faces are continuous. Such
functions, however, do not necessarily have continuous normal
component across the faces of the tetrahedrons. Thus, the
approximation $\vec{E}^h$ to $\vec{E}$ is not in $H(\dive)$ even
though $\vec{E}\in H(\dive)$.
 
\subsection{Discrete weak form}
After constructing the approximating spaces, the discrete problem is
constructed by restricting the bilinear form onto the piecewise
polynomial spaces.  In the following, we set $\vec{j}=0$ because we
are interested only in the dependence of the solution on initial
conditions.  Denoting
\[
H_h = H_{h,{\textrm{imp}}}(\curl)\times H_{h,{\textrm{imp}}}(\dive)\times H_{h,0}(\grad),
\]
and restricting~\eqref{mixed-form-1}--\eqref{mixed-form-3} to
$H_{h}$ leads to the following approximate variational problem:
Find $(\vec{E}^h,\vec{B}^h,p^h)\in H_h$ such that for all
$(\vec{F}^h,\vec{C}^h,q^h) \in H_{h}$ 
\begin{align}
&\scalare{\vec{E}_{t}^h}{\vec{F}^h} =-\scalare{\grad
  p^h}{\vec{F}^h}+\scalarmu{\vec{B}^h}{\curl\, \vec{F}^h}-(1+\gamma)\int_{\Gamma_{i}}
\scalar{\vec{n}\wedge \vec{E}^h}{\vec{n}\wedge \vec{F}^h} 
\label{mixed-form-11}\\
&\scalarmu{\vec{B}_{t}^h}{\vec{C}^h} =-\scalarmu{\curl\,
  \vec{E}^h}{\vec{C}},\label{mixed-form-21}\\
&(p_{t}^h, q^h) = \scalare{\vec{E}^h}{\grad q^h}.\label{mixed-form-31}
\end{align}
In the following, the superscript $h$ is omitted, since the considerations in the
rest of the paper are focused on the discrete problem in $H_h$. We also
interchangeably use $\vec{u}$ and $(\vec{E},\vec{B},\vec{p})^T$ and, similarly,
$\vec{w}$ and $(\vec{F},\vec{C},\vec{q})^T$. 

An operator is introduced, such that $\mathcal{A}:H_h\mapsto H_h$ via the bilinear
forms in~\eqref{mixed-form-11}--\eqref{mixed-form-31}.  For
$\vec{u}=(\vec{E},\vec{B},p)^T\in H_h$ and $\vec{w}=(\vec{F},\vec{C},q)^T\in H_h$, set
\begin{eqnarray*}
(\mathcal{A}\vec{u},\vec{w})
&=&-\scalare{\vec{E}}{\grad q}+\scalare{\grad p}{\vec{F}}
-\scalarmu{\vec{B}}{\curl\, \vec{F}}+
\scalarmu{\curl\, \vec{E}}{\vec{C}}.
\end{eqnarray*}
Corresponding to the boundary term, we also have the operator
associated with the impedance boundary condition, 
\[
(\mathcal{Z} \vec{u},\vec{w})
=(1+\gamma)\int_{\Gamma_{i}}\scalar{\vec{n}\wedge \vec{E}}{\vec{n}\wedge \vec{F}}
\]
Since we are now on a finite-dimensional space, we
write the semi-discrete problem (discretized in space and
continuous in time) as a constant coefficient linear system of ODEs,
i.e.
\begin{equation}\label{eqn:matsys}
\dot{\vec{u}}=-(\mathcal{A}+\mathcal{Z})
\vec{u}
\end{equation}
From the definitions of $\mathcal{A}$ and $\mathcal{Z}$ it is obvious
that $\mathcal{A}$ is skew symmetric and $\mathcal{Z}$ is symmetric
and positive semi-definite.  


\section{Matrix representation and time discretization}\label{sec:time}
We now show that the assembly of system \eqref{eqn:matsys} can be
constructed using only mass (Gramm) matrices formed with the bases in 
$H_{h,\textrm{imp}}(\curl)$, $H_{h,{\textrm{imp}}}(\dive)$, and $H_{h,0}(\grad)$ spaces and 
incidence matrices, whose entries encode the relationships
``vertex incident to an edge'', ``edge incident to a face'', etc.  The aim of this section is to provide some insight into the
implementation of such finite-element schemes and also to set the
stage for presenting the Crank--Nicolson discretization in time.

\subsection{Matrix representation}
We start with a description of the standard (canonical) bases in
$H_{h,0}(\grad)$, $H_{h,{\textrm{imp}}}(\dive)$, and $H_{h,{\textrm{imp}}}(\curl)$,
respectively.  By boundary vertices, edges, and faces we mean
vertices, edges, and faces lying on the boundary,
$\partial\Omega$. For an edge, this means that both its end vertices
are on the boundary and for a face it means that all three of its
vertices are on the boundary.  The remaining vertices, (edges, faces)
are designated as interior vertices (edges, faces).  We note that by a
standard convention, it is assumed that for the triangulation in hand
the directions of vectors tangential to edges and normal to faces are
fixed once and for all. It is easy and straightforward to check that a
change in these directions does not change the considerations that
follow.

We then have the following sets of degrees of freedom (DoFs):

\begin{itemize}
\item DoFs corresponding to the set of interior vertices
  $\{\vec{x}_i\}_{i=1}^{n_h}$:  A
  functional (also denoted by $\vec{x}_i$) is associated with an
  interior vertex $\vec{x}_i$ as
  $\vec{x}_i(q)=q(\vec{x}_i)$ for a sufficiently smooth function, $q$.

\item DoFs corresponding to the set of all interior edges
  and all edges on $\Gamma_i$: For a sufficiently smooth 
  vector-valued function, $\vec{v}$, and an edge,  $e\in \mathcal E$, the
  associated functional is $e(\vec{v})=\frac{1}{|e|}\int_e \vec{v}\cdot\vec{\tau}_e$, where
$\vec{\tau}_e$ is the unit vector tangential to the edge. The direction of
 the tangent vector, $\vec{\tau}_e$, is assumed to be fixed. 

\item DoFs corresponding to the set of interior faces,
  $\mathcal{F}$:  For a sufficiently smooth vector-valued
  function, $\vec{v}$, and a face $f\in \mathcal{F}$, the associated
  functional is $f(\vec{v}) = \frac{1}{|f|}\int_f \vec{v}\cdot
  \vec{n}_f$, where $\vec{n}_f$ is the unit vector normal to the face.
\end{itemize}

As bases for the spaces $H_{h,0}(\grad)$,
$H_{h,{\textrm{imp}}}(\curl)$, and $H_{h,{\textrm{imp}}}(\dive)$ we take the 
piecewise polynomial functions, which are dual to the functionals
given above. For the space $H_{h,0}(\grad)$, we denote these functions
by $\{\varphi_j\}_{j=1}^{n_h}$. They are piecewise linear,
continuous, and satisfy $\vec{x}_k(\varphi_j)=\delta_{kj}$, where
$\delta_{kj}$ is the Kroneker delta.

The bases for the other two spaces $H_{h,{\textrm{imp}}}(\curl)$ and
$H_{h,{\textrm{imp}}}(\dive)$ are then given in terms of the basis for
$H_{h,0}(\grad)$. For an edge $e\in \mathcal{E}$ with vertices
$(\vec{x}_i,\vec{x}_j)$ and a face $f\in \mathcal{F}$ with vertices
$(\vec{x}_i,\vec{x}_j,\vec{x}_{k})$,
\begin{eqnarray*}
\vec{\psi}_{e}&=&|e|(\varphi_i\grad\varphi_j-\varphi_j\grad\varphi_i),\\
\vec{\xi}_f&=&
|f|(\varphi_i(\grad\varphi_j\wedge\grad\varphi_k)
+\varphi_j(\grad\varphi_k\wedge\grad\varphi_i)
+\varphi_k(\grad\varphi_i\wedge\grad\varphi_j)).
\end{eqnarray*}
Here, $\vec{\tau}_e=(\vec{x}_j-\vec{x}_i)/|\vec{x}_i-\vec{x}_j|$ and 
the ordering of $(\vec{x}_i,\vec{x}_j,\vec{x}_{k})$ in a positive
direction is determined by the right-hand rule and the normal vector
$\vec{n}_f$.  We then have the following canonical representations of
functions in  
$H_{h,{\textrm{imp}}}(\curl)$, $H_{h,{\textrm{imp}}}(\dive)$, and 
$H_{h,0}(\grad)$:
\begin{eqnarray*}
\vec{v}\in  
H_{h,{\textrm{imp}}}(\curl),\quad & & 
\vec{v}=\sum_{e\in \mathcal{E}} e(\vec{v})
\vec{\psi}_{e}(\vec{x});\quad \vec{v}\in  H_{h,{\textrm{imp}}}(\dive),\quad 
\vec{v}=\sum_{f\in \mathcal{F}} f(\vec{v}) \vec{\xi}_f(\vec{x});\\
 q\in H_{h,0}(\grad),\quad &&
q=\sum_{i=1}^{n_h} \vec{x}_i(q) \varphi_i(\vec{x}).
\end{eqnarray*}
For functions that also depend on
time, i.e.  for the elements of $H_{h,{\textrm{imp}}}(\curl;t)$, $H_{h,{\textrm{imp}}}(\dive;t)$, and 
$H_{h,0}(\grad;t)$, we have similar representations with coefficients
depending on time as well. 
\begin{remark}
 In the rest of the paper, the same notation is used
  for the functions from $H_h$ and their vector representations in the
  bases given above. This is done in order to simplify the notation.
\end{remark}

The entries of the mass (Gramm) matrices for each of the piecewise polynomial spaces
are then,
\begin{equation*}
\left (\mathcal{M}_e\right )_{ee'} = 
(\vec{\psi}_{e},\vec{\psi}_{e^\prime}), \quad
\left (\mathcal{M}_f\right )_{ff'} =
(\vec{\xi}_f,\vec{\xi}_{f^\prime}),\quad
\left (\mathcal{M}_v\right )_{ij} = (\varphi_i,\varphi_j).
\end{equation*}

Next, the following matrix representations of the
operators defined in the previous section are introduced,
\begin{equation*}
\mathcal{G}_{ej}  = (\grad \varphi_j,\vec{\psi}_{e}),\quad
\mathcal{K}_{fe}= (\curl\, \vec{\psi}_{e},\vec{\xi}_f).
\end{equation*}
The matrix form of \eqref{eqn:matsys} is now rewritten as follows.
\begin{equation}\label{eqn:incidencesys}
\begin{small}
\left (
\begin{matrix}
\mathcal{M}_e&&\\
&\mathcal{M}_f&\\
&&\mathcal{M}_v\\
\end{matrix}\right) 
\left(\begin{matrix}
\dot{\vec{E}}\\
\dot{\vec{B}}\\
\dot{p}\\
\end{matrix} \right) = 
\left [ \left (
\begin{matrix}
&\mathcal{K}^T\mathcal{M}_f & -\mathcal{M}_e\mathcal{G}\\
-\mathcal{M}_f\mathcal{K}&&\\
\mathcal{G}^T\mathcal{M}_e&&\\
\end{matrix}\right) 
- \mathcal{Z} \right ]
\left(\begin{matrix}\vec{E}\\\vec{B}\\p\\\end{matrix}\right).
\end{small}
\end{equation}

\subsection{Time discretization}
To discretize system \eqref{eqn:incidencesys} in time, a
Crank--Nicolson scheme is used. We look at a time interval, $t\in [0,T]$, and
approximate the solution at $t=k\tau$, $k=1,\ldots,\frac{T}{\tau}$,
with $\tau$ a given time step.  Let $\vec{u}_k = \left (
  \vec{E}_k,\vec{B}_k,p_k \right)^T$ be the discrete approximation at
the current time $t=k\tau$, and $\vec{u}_{k-1} = \left (
  \vec{E}_{k-1},\vec{B}_{k-1},p_{k-1} \right)^T$ be the approximation
at the previous time $t=(k-1)\tau$.  Then, the Crank-Nicolson
formulation of \eqref{eqn:incidencesys} is
\begin{equation*}
\frac{1}{\tau}\mathcal{M} \left (
\vec{u}_k - \vec{u}_{k-1}\right ) = 
-\frac{1}{2} \left (\mathcal{A} + \mathcal{Z}) \right ) \left (
  \vec{u}_k + \vec{u}_{k-1} \right ),\quad \mbox{where}~\mathcal{M}=
\left (
\begin{matrix}
\mathcal{M}_e&&\\
&\mathcal{M}_f&\\
&&\mathcal{M}_v\\
\end{matrix}\right).
\end{equation*}
Rearranging the terms, we get the following linear system
for the approximate solution at time step $k\tau$ in terms of the
solution at time $(k-1)\tau$.
\begin{equation}\label{eqn:crank}
\left (\frac{1}{\tau} \mathcal{M} + \frac{1}{2}\left
  (\mathcal{A}+\mathcal{Z}\right ) \right )\vec{u}_k = \left (\frac{1}{\tau} \mathcal{M} - \frac{1}{2}\left
  (\mathcal{A}+\mathcal{Z}\right ) \right )\vec{u}_{k-1}.
\end{equation}

Next, we show that if the initial condition is weakly divergence-free,
as in the continuous case, we have that $p_k$ will remain zero for all
time and, thus, $\vec{E}_k$ is weakly divergence-free for all
$k$. This is the discrete analogue of
proposition~\ref{prop:weakly-div-free}.

\begin{lemma}\label{lemma:discrete-div-free}  
  Assume that $(\vec{E}_{0} , \grad q) = 0$ for all $q \in
  H_{h,0}(\grad)$.  For the Crank-Nicolson scheme described in
  \eqref{eqn:crank}, $p_{k}=0$ for all $k$ and $(\vec{E}_{k} , \grad q)
  = 0$ for all $q \in H_{h,0}(\grad)$ and all $k$. 
\end{lemma}

{\bf Proof.~~}
Start with $p_{0}=0$ and
$E_{0}$ being weakly divergence-free. It is shown that:
\[
\mbox{If} \quad p_{k} = 0 \quad \mbox{and} \quad
\mygrad^T\mass{e}\vec{E}_{k} = 0,
\quad\mbox{then}\quad 
p_{k+1} = 0 \quad \mbox{and} \quad
\mygrad^T\mass{e}\vec{E}_{k+1} = 0.
\]
This is the matrix representation of the assumptions and claims in the
lemma.  Setting $\alpha=2/\tau$ and using the defining relations for the
Crank-Nicolson time discretization, \eqref{eqn:incidencesys} and \eqref{eqn:crank}, the following linear system for
$\vec{E}_{k+1}$, $\vec{B}_{k+1}$, and $p_{k+1}$ is obtained:
\begin{align}
&\alpha\mass{e}\vec{E}_{k+1} - \mycurl^T\mass{f}\vec{B}_{k+1} + 
\mass{e}\mygrad p_{k+1} + \myz
\vec{E}_{k+1} \nonumber \\
&~~~~~~~~~~~~~~~~~~~ = \alpha\mass{e}\vec{E}_{k} +
\mycurl^T\mass{f}\vec{B}_{k} - \myz \vec{E}_{k}, \label{edges}\\
&\mass{f}\mycurl \vec{E}_{k+1} + \alpha\mass{f}\vec{B}_{k+1} = -\mass{f}\mycurl \vec{E}_{k} +
\alpha\mass{f}\vec{B}_{k}, \label{faces}\\
&-\mygrad^T\mass{e}\vec{E}_{k+1} + \alpha\mass{v}p_{k+1} = 0. \label{nodes}
\end{align}
Here, $\myz$ indicates the reduced matrix of $\mathcal{Z}$ applied to
just the edge DoF, $\vec{E}$.  Multiplying Equation \eqref{edges} from the left by $\mygrad^T$
yields,
\begin{eqnarray*}
\lefteqn{\alpha\mygrad^{T}\mass{e}\vec{E}_{k+1} -
\mygrad^{T}\mycurl^T\mass{f}\vec{B}_{k+1} 
+ \mygrad^{T}\mass{e}\mygrad p_{k+1} + \mygrad^{T}\myz
\vec{E}_{k+1}} \\
&&\qquad\qquad= \alpha\mygrad^{T}\mass{e}\vec{E}_{k} + 
\mygrad^{T}\mycurl^T\mass{f}\vec{B}_{k} - \mygrad^{T}\myz
\vec{E}_{k}.
\end{eqnarray*}
Next, note that $\mycurl\mygrad = 0$ (or equivalently
$\mygrad^T\mycurl^T=0$), since the curl of a gradient is zero.  Also,
since any $q \in H_{h,0}(\grad)\subset H_0^1(\Omega)$ is zero on
$\Gamma_i$ and $\Gamma_o$ and the tangential component of its gradient 
is zero on the boundary edges, then $\myz\mygrad = 0$ (or equivalently
$\mygrad^T\myz=0$, since $\myz$ is symmetric).  Thus, \eqref{edges}
simplifies to,
$$
\alpha\mygrad^{T}\mass{e}\vec{E}_{k+1} + \mygrad^{T}\mass{e}\mygrad
p_{k+1} = 0.
$$
Adding this to $\alpha$ times Equation
\eqref{nodes},  then gives
\begin{equation}\label{p1eqn}
\left ( \mygrad^{T}\mass{e}\mygrad  + \alpha^2\mass{v} \right) p_{k+1} = 0
\end{equation}
The above relation is the matrix representation of the variational problem:
\begin{equation*}
\left (\grad p_{k+1}, \grad q \right ) + \alpha^2\left ( p_{k+1},q \right) =
0,\quad \mbox{for all}\quad q\in H_{h,0}(\grad), 
\end{equation*}
and taking $q=p_{k+1}$ then gives that $p_{k+1} = 0$.  Finally, from
this fact and using \eqref{nodes}, it is immediately shown that
$\mygrad^T\mass{e}\vec{E}_{k+1} = 0$, concluding the proof.\hfill\qed

Thus, using the Crank-Nicolson scheme and appropriate initial
conditions, one can guarantee that the discrete approximation to the
electric field will be weakly divergence-free for all time.

\subsection{Solution of the discrete linear systems} 
To solve the system, we look at the matrix corresponding to $ \frac{1}{\tau}
\mathcal{M} + \frac{1}{2}\left (\mathcal{A}+\mathcal{Z}\right )$, which
is on the left side of~\eqref{eqn:crank}. We have to solve the system
with this matrix at every time step.    
Using the incidence matrices as in
\eqref{eqn:matsys}, this operator is written as
\begin{equation*}
\frac{1}{\tau} \mathcal{M} + \frac{1}{2}\left
  (\mathcal{A}+\mathcal{Z}\right ) = 
\frac12\left (\begin{matrix}
\frac{2}{\tau}\mathcal{M}_e&-\mathcal{K}^T\mathcal{M}_f & \mathcal{M}_e\mathcal{G}\\
\mathcal{M}_f\mathcal{K}&\frac{2}{\tau}\mathcal{M}_f&\\
-\mathcal{G}^T\mathcal{M}_e&&\frac{2}{\tau}\mathcal{M}_v\\
\end{matrix}\right )
+ \frac{1}{2}\mathcal{Z} .
\end{equation*}
Since, the mass matrices, $\mathcal{M}_e$, $\mathcal{M}_f$,
$\mathcal{M}_v$, are all SPD and $\mathcal{Z}$ is symmetric positive
semi-definite and only contributes to the edge-edge diagonal block of
the system, the entire operator can be made symmetric by a simple
permutation.  Multiplying on the left by $J = \left
  ( \begin{matrix}I&&\\&-I&\\&&-I\\\end{matrix}\right)$ will yield the
operator
\begin{equation*}
J\left (\frac{1}{\tau} \mathcal{M} + \frac{1}{2}\left
  (\mathcal{A}+\mathcal{Z}\right ) \right ) = 
\left (\begin{matrix}
\frac{1}{\tau}\mathcal{M}_e&-\frac{1}{2}\mathcal{K}^T\mathcal{M}_f & +\frac{1}{2}\mathcal{M}_e\mathcal{G}\\
-\frac{1}{2}\mathcal{M}_f\mathcal{K}&-\frac{1}{\tau}\mathcal{M}_f
  &\\
\frac{1}{2}\mathcal{G}^T\mathcal{M}_e&&-\frac{1}{\tau}\mathcal{M}_v\\
\end{matrix}\right )
+ \frac{1}{2}\mathcal{Z} ,
\end{equation*}
which is now symmetric.  Therefore, the final system to solve is
\begin{equation}\label{eqn:solvesys}
\mathcal{J}\left (\frac{1}{\tau} \mathcal{M} + \frac{1}{2}\left
  (\mathcal{A}+\mathcal{Z}\right ) \right )\vec{u}_n = \mathcal{J}\left (\frac{1}{\tau} \mathcal{M} - \frac{1}{2}\left
  (\mathcal{A}+\mathcal{Z}\right ) \right )\vec{u}_{n-1},
\end{equation}
and a standard iterative solvers such as MINRES can be applied.  This is
used in the test problems below.  

\section{Numerical Results}\label{sec:numerics}
Here, we perform some numerical tests by solving system
\eqref{eqn:solvesys} using the Crank-Nicolson time discretization
described in the previous section.  To test for decay in the energy of
the solution, we start with initial conditions and boundary conditions
of the form described in \cite{2011ColombiniF_PetkovV_RauchJ-aa}.  We
take as the domain, the area between a polyhedral approximation of the
sphere of radius 1 and a polyhedral approximation of a sphere of
radius 4.  The inner sphere represents the impedance boundary and the
outer sphere is considered far enough away (and it is for the solutions we
approximate) that a Dirichlet-like perfect conductor boundary
condition on the outer sphere is used.  In other words, we prescribe $\vec{E}\wedge n$, $\vec{B}\cdot n$, and
$p=0$ on the outer sphere.  The exact solution (taken
from~\cite[Theorem~3.2]{2011ColombiniF_PetkovV_RauchJ-aa}) is given as
follows:
\begin{eqnarray}
\vec{E}_* & =& \frac{e^{r\left (|\vec{x}| + t\right )}}{|\vec{x}|^2}\left ( r^2 -
  \frac{r}{|\vec{x}|}\right ) \left (\begin{array}{c}0\\z\\-y\\\end{array}
\right ),\label{Einit}\\
\vec{B}_* &=& e^{r\left (|\vec{x}| + t\right )} \left [ \frac{1}{|\vec{x}|^3}\left ( r^2 -
  \frac{3r}{|\vec{x}|} + \frac{3}{|\vec{x}|
^2}\right ) \left (\begin{array}{c}z^2+y^2\\-xy\\-xz\\\end{array}
\right ) +
\left ( \begin{array}{c}\frac{2r}{|\vec{x}|} -
    \frac{2}{|\vec{x}|^2}\\0\\0\\\end{array} \right ) \right ].\label{Binit}
\end{eqnarray}
Different values of $\gamma$ yield different values of $r$ in
solutions
\eqref{Einit}-\eqref{Binit}. Following~\cite{2011ColombiniF_PetkovV_RauchJ-aa},
we have that $(\vec{E}_*,\vec{B}_*,0)$ solves Maxwell's system
with an impedance boundary condition on $\Gamma_i$ and
$$r = 1/2\left ( 1-\sqrt{1+4/\gamma}\right ).$$
For the tests below, we take $\gamma = 0.05$ ($r = -4$).  

For the annular domain that we consider, a basis in the one
  dimensional space of harmonic forms is of the form
  $\grad\hh=\vec{x}/|\vec{x}|^3$. Direct computation shows that 
$\langle\vec{E}_*,\grad \hh\rangle=0$, and this identity holds pointwise and for
  all $t$. Thus, the ADS that we are trying to approximate is
  orthogonal to the harmonic forms, in this particular example.

\subsection{Approximation of the initial conditions} 
Since the solutions given above are not in the finite-dimensional
spaces considered, we take an initial condition $\vec{E}_{0}$, which
is based on the piecewise polynomial interpolant of the
exponentially-decaying solution given in equations~\eqref{Einit} at
$t=0$.  In other words, we set
\begin{equation}\label{E0init}
\widetilde{\vec{E}}_{0} = 
\sum_{e\in \mathcal{E}}
e(\vec{E}_*(0,\vec{x}))\vec{\psi}_e(\vec{x}). 
\end{equation}

  \begin{remark} We noted earlier that the ADS given by~\eqref{Einit}
    is orthogonal to all harmonic forms. However, this is not so for
    the discrete approximation in~\eqref{E0init}. We perform an
    orthogonalization step and take initial conditions that are
    orthogonal to the \emph{discrete} harmonic forms in order to
    approximate the ADS  and the energy decay more accurately. 
\end{remark}

 We further correct $\widetilde{\vec{E}}_{0}$ to get an initial guess
that is orthogonal to the gradients as well as the gradients of the
discrete harmonic form, $\grad \hh^h$. This is done in a standard
fashion by projecting out these gradients. First, we find $s\in
H_{0,h}(\grad)$, such that for all $q\in H_{0,h}(\grad)$, we have
$$(\grad s,\grad q) = (\widetilde{\vec{E}}_0,\grad q),\quad 
\vec{E}_{0} = \widetilde{\vec{E}}_{0,h} - 
\grad s- \frac{(\widetilde{\vec{E}}_{0,h},\grad\hh^h)}{\|\grad\hh^h\|^2}\grad\hh^h.
$$
As a result, $\vec{E}_{0}$ is orthogonal to the gradients of functions
in~$H_{0,h}(\grad)$ and also to the gradient of the discrete harmonic
form. We note that this orthogonalization requires two solutions of
Laplace equation.  It is straightforward to see that if the initial
guess satisfies such condition, then the solution $\vec{E}$ satisfies
this condition for all times.  Finally, $\vec{B}_0$ is computed as
$\vec{B}_0 = \frac{1}{r} \mathcal{K}\vec{E}_{0}$. 


\subsection{Numerical results}
We test the approximation to the asymptotically disappearing solutions
on  a grid with 728 ($h=1/8$), 4,886 ($h=1/16$), 35,594 ($h=1/32$), and 271,250 ($h=1/64$) 
nodes on the domain.  The computational domain is shown in
Figure \ref{fig:sphere}.  We run a MINRES solver on the Crank-Nicolson system that is
symmetrized, Equation \eqref{eqn:solvesys}, for 20 time
steps using a step size $\tau=0.1$.  

The results are shown in
Tables~\ref{table:sphere:h8}--\ref{table:sphere:h64}, where we display
the $\|\vec{E}\|_{L_2(\Omega)}$ and $\|\vec{B}\|_{L_2(\Omega)}$ norms. The
total energy of this system is given by
$\|\vec{E}\|^2_{L_2(\Omega)}+\|\vec{B}\|^2_{L_2(\Omega)}$. 
\begin{table}[h!]
\begin{minipage}[h!]{0.45\linewidth}
\centering
{\small
\begin{tabular}{|c|cc|}
\hline 
Step&$\|\vec{E}\|_{L_2(\Omega)}$&$\|\vec{B}\|_{L_2(\Omega)}$\\
\hline
0&0.906&0.559\\
1&0.830&0.334\\
2&0.723&0.155\\
3&0.614&0.116\\
4&0.526&0.172\\
5&0.450&0.242\\
6&0.380&0.304\\
7&0.316&0.349\\
8&0.262&0.379\\
9&0.233&0.389\\
10&0.234&0.384\\
11&0.248&0.373\\
12&0.266&0.358\\
13&0.284&0.339\\
14&0.292&0.322\\
15&0.288&0.312\\
16&0.276&0.307\\
17&0.262&0.303\\
18&0.251&0.299\\
19&0.248&0.293\\
20&0.246&0.287\\
\hline
\end{tabular}
}
\caption{Sphere Obstacle.  $\gamma =
  0.05$. $h = 1/8$.}
\label{table:sphere:h8}
\end{minipage}
\hspace{0.5cm}
\begin{minipage}[h!]{0.45\linewidth}
\centering
{\small
\begin{tabular}{|c|cc|}
\hline
Step&$\|\vec{E}\|_{L_2(\Omega)}$&$\|\vec{B}\|_{L_2(\Omega)}$\\
\hline
0&0.553&0.426\\
1&0.451&0.266\\
2&0.349&0.153\\
3&0.267&0.093\\
4&0.200&0.086\\
5&0.147&0.103\\
6&0.117&0.103\\
7&0.107&0.099\\
8&0.104&0.097\\
9&0.095&0.101\\
10&0.088&0.103\\
11&0.085&0.104\\
12&0.084&0.104\\
13&0.085&0.102\\
14&0.086&0.100\\
15&0.086&0.099\\
16&0.086&0.099\\
17&0.085&0.099\\
18&0.086&0.097\\
19&0.085&0.095\\
20&0.085&0.095\\
\hline
\end{tabular}
}
\caption{Sphere Obstacle.  $\gamma =
  0.05$. $h = 1/16$.}
\label{table:sphere:h16}
\end{minipage}
\end{table}

\begin{table}[h!]
\begin{minipage}[h!]{0.45\linewidth}
\centering
{\small
\begin{tabular}{|c|cc|}
\hline 
Step&$\|\vec{E}\|_{L_2(\Omega)}$&$\|\vec{B}\|_{L_2(\Omega)}$\\
\hline
0&0.425&0.405\\
1&0.308&0.262\\
2&0.214&0.163\\
3&0.148&0.102\\
4&0.103&0.061\\
5&0.074&0.044\\
6&0.057&0.039\\
7&0.046&0.037\\
8&0.040&0.035\\
9&0.035&0.034\\
10&0.033&0.034\\
11&0.032&0.034\\
12&0.031&0.035\\
13&0.031&0.034\\
14&0.031&0.034\\
15&0.031&0.034\\
16&0.031&0.034\\
17&0.031&0.033\\
18&0.031&0.033\\
19&0.030&0.033\\
20&0.031&0.032\\
\hline
\end{tabular}
}
\caption{Sphere Obstacle.  $\gamma =
  0.05$. $h = 1/32$.}
 \label{table:sphere:h32}
\end{minipage}
\hspace{0.5cm}
\begin{minipage}[h!]{0.45\linewidth}
\centering
{\small
\begin{tabular}{|c|cc|}
\hline
Step&$\|\vec{E}\|_{L_2(\Omega)}$&$\|\vec{B}\|_{L_2(\Omega)}$\\
\hline
0&0.383&0.401\\
1&0.261&0.264\\
2&0.176&0.172\\
3&0.120&0.113\\
4&0.081&0.075\\
5&0.056&0.051\\
6&0.039&0.034\\
7&0.028&0.025\\
8&0.021&0.020\\
9&0.017&0.017\\
10&0.015&0.015\\
11&0.014&0.015\\
12&0.013&0.015\\
13&0.012&0.015\\
14&0.013&0.014\\
15&0.013&0.014\\
16&0.013&0.014\\
17&0.012&0.015\\
18&0.013&0.014\\
19&0.013&0.014\\
20&0.013&0.014\\
\hline
\end{tabular}
}
\caption{Sphere Obstacle.  $\gamma =
  0.05$. $h = 1/64$.}
\label{table:sphere:h64}
\end{minipage}
\end{table}

\begin{figure}[h!]
\centering
\includegraphics[scale=0.35]{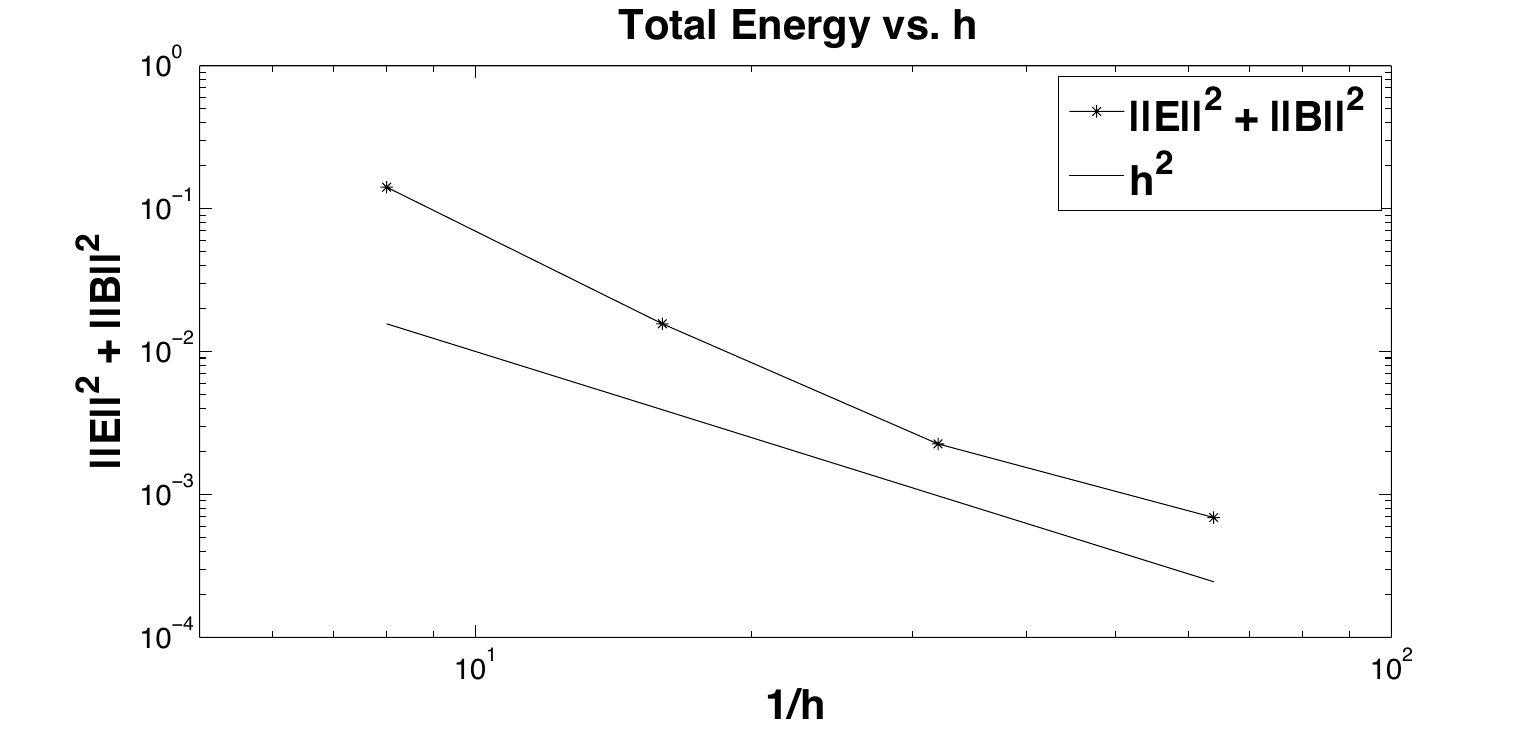}
\label{fig:maxenergy}
\caption{Plot of Total Energy
  ($\|E\|^2_{L_2(\Omega)}+\|B\|^2_{L_2(\Omega)}$) vs. mesh size after
  20 time steps.  The
  x-axis shows $1/h$.}
\end{figure}

Tables \ref{table:sphere:h8}-\ref{table:sphere:h64} show that over
time the $L_2$ norms of the electric and magnetic fields decay.  For
each mesh size, it appears that the energy reaches some steady-state
value, where it does not decay anymore.  Figure \ref{fig:maxenergy}
shows that this final energy value (at time step 20),
$||E||^2_{L_2(\Omega)}+||B||^2_{L_2(\Omega)}$, decreases as $h^2$,
when $h\to 0$.  Thus, as the polyhedron domain more closely represents
the spherical domain, as in \cite{2011ColombiniF_PetkovV_RauchJ-aa},
the total energy should decay to zero as expected over time.

\section{Concluding remarks}\label{sec:conclusion}
We have shown that using appropriate finite-element spaces one can
approximate the asymptotically-decaying solutions to Maxwell's
equations on the exterior of a spherical obstacle.
We also have shown that the Crank-Nicolson time discretization keeps
the electric and magnetic fields divergence free (the electric field
weakly and the magnetic field--strongly).  The next step is to apply
the methods to more general domains and answer the question of whether
ADS exist for obstacles with more complicated geometry.  In relation
to considering less regular obstacles, there are many open problems
for the analysis of systems with dissipative boundary conditions.  The
computational techniques that we have introduced here easily
generalize to cases of variable, matrix valued permittivity and permeability
$\varepsilon(x)$ and $\mu(x)$, as well as to more general hyperbolic
systems for which the ADS phenomenon occurs.

\section*{Acknowledgments}
We are grateful to anonymous referees which provided many valuable
remarks and helped us improve the exposition. 
This work was supported in part by the U.S. Department of Energy grant
DE-FG02-11ER26062/DE-SC0006903 and subcontract LLNL-B595949. The
third author was supported in part by the National Science
Foundation grants DMS-0810982 and DMS-1217142.

\bibliographystyle{siam}

\end{document}